\documentclass[a4paper,12pt]{amsart}

\usepackage{graphicx}

\newtheorem{theorem}{Theorem}[section]
\newtheorem{proposition}[theorem]{Proposition}

\newtheorem{lemma}[theorem]{Lemma}
\newtheorem{conjecture}[theorem]{Conjecture}
\newtheorem{metatheorem}[theorem]{Metatheorem}

\newtheorem{Example}[theorem]{Example}
\newenvironment{example}{\begin{Example}\rm}{\end{Example}}
\newtheorem{Examples}[theorem]{Examples}

\newtheorem{Remark}[theorem]{Remark}
\newenvironment{remark}{\begin{Remark}\rm}{\end{Remark}}
\newtheorem{Remarks}[theorem]{Remarks}

\newtheorem{Question}[theorem]{Question}

\newtheorem{Summary}[theorem]{Summary}
\newenvironment{summary}{\begin{Summary}\rm}{\end{Summary}}

\newcommand{\al}{\alpha}
\newcommand{\be}{\beta}
\newcommand{\ga}{\gamma}
\newcommand{\Ga}{\Gamma}

\newcommand{\la}{\lambda}
\newcommand{\La}{\Lambda}

\newcommand{\si}{\sigma}
\newcommand{\Si}{\Sigma}
\newcommand{\ze}{\zeta}

\let\cal=\mathcal
\let\Bbb=\mathbb

\begin{document}

\title[Prime-pair constants, 10.06.08]{Mean value one of 
prime-pair constants}


\date{June 10, 2008}

\author{Fokko van de Bult and Jaap Korevaar}

\begin{abstract}
For $k>1$, $r\ne 0$ and large $x$, let $\pi^k_{2r}(x)$
denote the number of prime pairs $(p,\,p^k+2r)$ with
$p\le x$. By the Bateman--Horn conjecture the function
$\pi^k_{2r}(x)$ should be asymptotic to
$(2/k)C^k_{2r}{\rm li}_2(x)$, with certain specific
constants $C^k_{2r}$. Heuristic arguments lead to the
conjecture that these constants have mean value one, just
like the Hardy--Littlewood constants $C_{2r}$ for prime
pairs $(p,\,p+2r)$. The conjecture is supported by
extensive numerical work.
\end{abstract}

\maketitle

\setcounter{equation}{0}    
\section{Introduction} \label{sec:1}
In the following $p$ runs through the primes. We are
interested in the `Bateman--Horn constants' associated
with prime pairs $(p,\,p^k+2r)$, where
$k\ge 2$ and $r$ runs over ${\Bbb Z}\setminus 0$. 
  
For convenience we first state the general Bateman--Horn
conjecture. It involves an $m$-tuple
$f=\{f_1,\cdots,f_m\}$ of polynomials $f_j$ with integer
coefficients and nonconstant ratios, and of positive
degrees $d_1,\cdots,d_m$. The conjecture involves the
counting function
\begin{equation} \label{eq:1.1}
\pi_f(x)=\#\{1\le n\le x:\,f_1(n),\cdots,f_m(n)
\;\mbox{all prime}\}.
\end{equation} 
Let
\begin{equation} \label{eq:1.2}
N_f(p)=\#\{n,\,1\le n\le p:\,f_1(n)\cdots f_m(n)\equiv
0\; ({\rm mod}\,p)\}.
\end{equation}
Assuming that the polynomials $f_j(n)$ are
irreducible and that $N_f(p)<p$ for every prime $p$,
Schinzel and Sierpinski \cite{ScS58} had conjectured that
$\pi_f(x)\to\infty$ as $x\to\infty$. The corresponding
quantitative conjecture is due to Bateman and Horn
\cite{BH62}, \cite{BH65}. Without making the above
assumptions, let
\begin{equation} \label{eq:1.3}
BH(f)=\frac{1}{d_1\cdots
d_m}\prod_p\,\left(1-\frac{1}{p}\right)^{-m}
\left(1-\frac{N_f(p)}{p}\right).
\end{equation}
The product will converge, but $BH(f)$ may be zero;
if one of the polynomials $f_j(n)$ can be factored, one
may define $BH(f)=0$. The conjecture now reads as follows;
cf.\  also Schinzel \cite{Sc63}, Davenport and Schinzel
\cite{DS66}, and the recent survey paper by Hindry and
Rivoal \cite{HiRi05}.
\begin{conjecture} \label{con:1.1}
Let $\pi_f(x)$ be as in $(\ref{eq:1.1})$. Then 
\begin{equation} \label{eq:1.4}
\pi_f(x)\sim BH(f)\,{\rm
li}_m(x)=BH(f)\,\int_2^x\frac{dt}{\log^m t},
\end{equation}
in the sense that $\pi_f(x)/{\rm li}_m(x)\to BH(f)$
as $x\to\infty$.
\end{conjecture}
One may verify that for prime pairs
$(p,\,p+2r)$ (with $r\in{\Bbb N}$), this
gives the classical conjecture of Hardy and Littlewood
\cite{HL23}:
\begin{equation} \label{eq:1.5}
\pi_{2r}(x)=\#\{p\le
x:\,p+2r\;\mbox{prime}\}\sim 2C_{2r}{\rm li}_2(x),
\end{equation}
where
\begin{equation} \label{eq:1.6}
C_2 = \prod_{p\,{\rm
prime},\,p>2}\,\left\{1-\frac{1}{(p-1)^2}\right\},\quad
C_{2r} = C_2\prod_{p|r,\,p>2}\frac{p-1}{p-2}.
\end{equation}

Turning to prime pairs $(p,\,p^k+2r)$, where
$k\ge 2$ and $r\in{\Bbb Z}\setminus 0$, we consider the
pair of polynomials
\begin{equation} \label{eq:1.7} 
f_{2r}(n)=f^k_{2r}(n)=\{n,\,n^k+2r\}.
\end{equation}
Adjusting our earlier notation, it is convenient to write
\begin{align} \label{eq:1.8}
\pi_{2r}(x) &= \pi^k_{2r}(x) = \pi_{f_{2r}}(x) = \#\{p\le
x:\,p^k+2r\;\mbox{prime}\}, \notag \\
N_{2r}(p) &= N^k_{2r}(p) =\#\{n,\,1\le
n\le p:\,n(n^k+2r)\equiv 0\; ({\rm mod}\,p)\}.
\end{align}
Suitable constants are defined by
\begin{equation} \label{eq:1.9}
C_{2r}=C^k_{2r} =
\prod_{p>2}\,\left(\frac{p}{p-1}\right)^2
\frac{p-N_{2r}(p)}{p};
\end{equation}
we set $C^k_{2r}=0$ if $n^k+2r$ can be factored. The
corresponding BH constants would be $(2/k)C^k_{2r}$. It
is known that the Hardy--Littlewood constants
$C_{2r}=C^1_{2r}$ for prime pairs $(p,\,p+2r)$ have
average one; cf.\ Section \ref{sec:4}. Our paper
provides both heuristic and numerical support for an
extension involving prime pairs $(p,\,p^k+2r)$:
\begin{metatheorem} \label{the:1.2}
For any degree $k\ge 2$, the adjusted Bateman--Horn
constants $C^k_{2r}$ have mean value one:
\begin{equation} \label{eq:1.10}
S^k_\la=\sum_{0<|2r|\le\la}\, C^k_{2r}
\sim \la\;\;\mbox{as}\;\; \la\to\infty.
\end{equation}
\end{metatheorem}
It is convenient to introduce the following auxiliary
functions $g_q(n)$, where $q\in{\Bbb Z}\setminus 0$ will
usually be taken equal to $2r$:
\begin{align} & \label{eq:1.11} 
g_q(n) = g^k_q(n)=n^k+q, \notag \\ &
\nu_q(p) = \nu^k_q(p)=\#\{n,\,1\le
n\le p:\,g_q(n)\equiv 0\; ({\rm mod}\,p)\}, \\ &
\ga_q = \ga^k_q=
\prod_{p>2}\,\frac{p}{p-1}\,\frac{p-\nu_q(p)}{p}.\notag
\end{align}
If $g_q(n)$ can be factored we set $\ga_q=0$.
Observe that $N_{2r}(p)=\nu_{2r}(p)+1$ except when
$p|2r$; in the latter case $N_{2r}(p)=\nu_{2r}(p)$.
Thus if $\ga_{2r}\ne 0$, the ratio $C_{2r}/\ga_{2r}$ is
given by an absolutely convergent product:
\begin{equation} \label{eq:1.12}
C_{2r}/\ga_{2r}=\prod_{p|2r,\,p>2}\,\frac{p}{p-1}
\prod_{p\not|\,2r}\,\frac{p}{p-1}\,
\frac{p-\nu_{2r}(p)-1}{p-\nu_{2r}(p)}.
\end{equation}

\setcounter{equation}{0} 
\section{Prime pairs $(p,\,p^2+2r)$}
\label{sec:2}
Consider the pair of polynomials
\begin{equation} \label{eq:2.1}
f(n)=f_{2r}(n)=\{n,\,g_{2r}(n)\}=\{n,\,n^2+2r\}\qquad
(r\in{\Bbb Z}\setminus 0).
\end{equation}
The functions $\pi^2_{2r}(x)=\pi_{f_{2r}}(x)$ will remain
bounded as $x\to\infty$ if $g_{2r}(n)$ can be factored, or
if $r\equiv 1$ (mod $3$); in the latter case $n(n^2+2r)$
is always divisible by $3$. For primes $p$ the
number of solutions of the quadratic congruence
$$n^2\equiv -2r\;({\rm mod}\;p)$$
is given by
$$\nu_{2r}(p)=1+\left(\frac{-2r}{p}\right),$$
where $(-2r/p)$ is the Legendre symbol. It follows that
\begin{equation} \label{eq:2.2}
N_{2r}(p)=2+\left(\frac{-2r}{p}\right)
\quad\mbox{for}\;\;p\not|\,2r,
\end{equation}
while $N_{2r}(p)=\nu_{2r}(p)=1$ if $p|2r$. The values
$\chi(p)=(-2r/p)$ generate a real Dirichlet character
(different from the principal character) belonging to
some modulus
$m=m_{2r}$. The convergence of the products for
$\ga^2_{2r}$ and $C^2_{2r}$ now follows from the
convergence of the series 
$$\sum_{p>2}\,\frac{\nu_{2r}(p)-1}{p}
=\sum_{p>2}\,\frac{\chi(p)}{p},$$
which is a classical result for Dirichlet characters; cf.\
Landau \cite{La09}.

\medskip\noindent
{\scshape The special case} $2r=-2$ or $g(n)=n^2-2$.
For $p>2$
$$\nu_{-2}(p) = 1+\left(\frac{2}{p}\right)
=\left\{\begin{array}{ll} 2 &
\mbox{if $p\equiv\pm 1$\;({\rm mod}\;8)},
\\ 0 & \mbox{otherwise}.
\end{array}\right.$$
In this case the values
$$\chi(2)=0\;\;\mbox{and}\;\;\chi(p)=(-1)^{(p^2-1)/8}
\;\;\mbox{for}\;\;p>2$$
generate a character modulo $8$.
\begin{table} \label{table:1}
\begin{tabular}{lllll} 
$x$     & $\pi^2_{-2}(x)$ & $L_2(x)$    & $\rho(x)$ & \\
        &                 &             &        & \\
$10$    &  4              &             &        & \\
$10^2$  & 13              &             &        & \\
$10^3$  & 52              &             &        & \\
$10^4$  & 259             & 274         & 0.945  & \\ 
$10^5$  & 1595            & 1600        & 0.997  & \\
$10^6$  & 10548           & 10567       & 0.998  & \\
$10^7$  & 74914           & 75275       & 0.995  & \\
$10^8$  & 563533          & 564200      & 0.999  & \\
        &                 &             &        &       
\end{tabular}
\caption{Counting prime pairs $(p,\,p^2-2)$}
\end{table}  

A rough computation shows that here
\begin{equation} \label{eq:2.3}
BH(f)=C^2_{-2}\approx 1.6916.
\end{equation}
We have also computed the counting function 
\begin{equation} \label{eq:2.4}
\pi^2_{-2}(x)=\#\{p\le x:\,p^2-2\;\mbox{prime}\}
\end{equation} 
for $x=10,\,10^2,\,\cdots,\,10^8$. In
Table $1$ the number $\pi^2_{-2}(x)$ is compared to
rounded values 
$$L_2(x)\;\;\mbox{of}\;\;1.6916\;{\rm
li}_2(x)=1.6916\int_2^x\,\frac{dt}{\log^2 t}.$$ 
The table includes some ratios 
$$\rho(x)=\pi^2_{-2}(x)/L_2(x).$$
These seem to converge to $1$ rather quickly!
\begin{table} \label{table:2}
\begin{tabular}{llllll} 
$2r$\hspace{1cm} & $\ga^2_{2r}$\hspace{1cm} &
$C^2_{2r}$\hspace{1cm} & $\ga^2_{-2r}$\hspace{1cm} &
$C^2_{-2r}$   & \\
       &        &         &        &           & \\
$2$    & 0.71   & 0       & 1.85   & 1.692     & \\
$4$    & 1.37   & 1.107   & 0      & 0         & \\
$6$    & 0.71   & 0.806   & 1.04   & 1.270     & \\
$8$    & 0.71   & 0       & 1.85   & 1.692     & \\ 
$10$   & 1.08   & 1.194   & 0.67   & 0         & \\
$12$   & 1.12   & 1.522   & 1.38   & 1.976     & \\
$14$   & 0.42   & 0       & 1.15   & 1.070     & \\
$16$   & 1.37   & 1.107   & 0      & 0         & \\
$18$   & 1.43   & 2.048   & 1.23   & 1.692     & \\
$20$   & 0.53   & 0       & 1.77   & 2.131     & \\ 
$22$   & 1.77   & 1.872   & 0.60   & 0         & \\
$24$   & 0.71   & 0.806   & 1.04   & 1.270     & \\
$26$   & 0.37   & 0       & 1.17   & 1.007     & \\
$28$   & 1.97   & 2.220   & 0.78   & 0         & \\
$30$   & 0.87   & 1.532   & 0.86   & 1.450     & \\
       &        &         &        &           &      
\end{tabular}
\caption{Bateman--Horn constants for $k=2$}
\end{table}

\setcounter{equation}{0} 
\section{Prime pairs $(p,\,p^3\pm 2r)$}
\label{sec:3}
For $r\in{\Bbb N}$ we now consider the pairs of
polynomials
\begin{equation} \label{eq:3.1}
f(n)=f_{2r}(n)=\{n,\,g_{2r}(n)\}=\{n,\,n^3\pm 2r\}.
\end{equation}
Modulo $p$ the number $\nu(p)$ of solutions of the cubic
congruence
\begin{equation} \label{eq:3.2}
n^3\equiv q\;({\rm mod}\;p)\qquad(q\in{\Bbb Z\setminus 0})
\end{equation}
is equal to $1,\,0$ or $3$, depending on $|q|$; cf.\ 
Ireland and Rosen \cite{IR90}.
If $p|q$ one has $\nu(p)=1$; in the case $q=2r$ and $p|2r$
also $N(p)=1$. Ignoring such $p$ for the moment, the
primes divide into three classes. The class of
`$1$-primes' is independent of $q$. It consists of $p=3$
and the primes $p\equiv 2$ (mod $3$), cf.\ Sloane
\cite{Sl01}:
\begin{align} &
2,\,3,\,5,\,11,\,17,\,23,\,29,\,41,\,47,\,53,\,59,\,71,
\,83,\,89,\notag \\ &
101,\,107,\,113,\,131,\,137,\,149,\,167,\,173,\,179,\,191,
\,197, \notag \\ &
227,\,233,\,239,\,251,\,257,\,263,\,269,\,281,\,293,\,311,
\,\cdots.\notag 
\end{align}
For these $p$, roughly half of the primes, equation
(\ref{eq:3.2}) always has precisely one solution $n$. 

The primes $p\equiv 1$ (mod $3$) are `unstable'. The
convergence of the product for $\ga^3_q$ in
(\ref{eq:1.11}) shows that $\nu(p)=0$ for roughly two
thirds of these primes and $\nu(p)=3$ for roughly
one third of them. The classes of `$0$-primes' and
`$3$-primes' depend on $|q|$. 
\begin{example} \label{exam:3.1} {\scshape The special
case} $f(n)=\{n,\,n^3\pm 2\}$. The $3$-primes in the case
$2r=\pm 2$ were characterized by Euler and Gauss; cf.\
Cox \cite{Cox89}. They are the primes of the form
$p=a^2+27b^2$, cf.\ Sloane \cite{Sl06} and additional
references in Section \ref{sec:4}: 
$$31,\,43,\,109,\,127,\,157,\,223,\,229,\,277,\,283,\,
307,\,\cdots.$$
The remaining primes $p\equiv 1$ (mod $3$) are
$0$-primes, cf.\ Sloane \cite{Sl05}:
\begin{align} & 
7,\,13,\,19,\,37,\,61,\,67,\,73,\,79,\,97,\,103,\,139,
\notag \\ &
151,\,163,\,181,\,193,\,199,\,211,\,241,\,271,\,313,\,
\cdots. \notag
\end{align}
Using the primes $p$ up to large $N$, formulas
(\ref{eq:1.11}) and (\ref{eq:1.12}) give
\begin{align} 
\ga^3_2 &\approx \prod_{p\le N,\,\nu(p)=0}\,\frac{p}{p-1}
\prod_{p\le N,\,\nu(p)=3}\,\frac{p-3}{p-1}\approx 1.30,
\notag\\ C^3_2 &=
\ga^3_2\,\prod_{\nu(p)=1;\,p>2}\,
\frac{p}{p-1}\,\frac{p-2}{p-1}\,
\prod_{\nu(p)=3}\,\frac{p}{p-1}\,\frac{p-4}{p-3}
\approx 0.89. \notag
\end{align}
\end{example}
\begin{example} \label{exam:3.2} {\scshape The special
case} $f(n)=\{n,\,n^3\pm 10\}$. Using the $3$-primes
$p<500$ from Table $3$ and paying special attention to
$p=5$ one finds
\begin{align} \label{eq:3.4}
\ga^3_{10} &\approx
\prod_{p<500,\,\nu(p)=0}\,\frac{p}{p-1}
\prod_{p<500,\,\nu(p)=3}\,\frac{p-3}{p-1}\approx 1.34,
\notag\\ C^3_{10} &=
\ga^3_{10}\cdot\frac{5}{4}\,\prod_{\nu(p)=1;\,p\ne
2,\,5}\,\frac{p}{p-1}\,\frac{p-2}{p-1}\,
\prod_{\nu(p)=3}\,\frac{p}{p-1}\,\frac{p-4}{p-3}
\approx 1.22.\notag
\end{align}
\end{example}

\setcounter{equation}{0} 
\section{Average of constants $C^k_{2r}$}
\label{sec:4}
The constants $C_{2r}=C^1_{2r}$ associated with
ordinary prime-pairs $(p,\,p+2r)$, with $r\in{\Bbb
N}$, have mean value one:
\begin{proposition} \label{prop:4.1}
One has
\begin{equation} \label{eq:4.1}
S_m=\sum_{1\le r\le m} C_{2r}\sim m\;\;\mbox{as}\;\;
m\to\infty.
\end{equation}
\end{proposition}
An extension to the constants in the `prime $n$-tuple
conjecture' was given by Gallagher \cite{Ga76}. Strong
estimates for the sums $S_m$ are due to  Bombieri and
Davenport \cite{BD66}, Montgomery
\cite{Mo71}, and Friedlander and Goldston \cite{FG95}.
Using singular series the latter showed that
\begin{equation} \label{eq:4.2}
S_m=m-(1/2)\log m+\cal{O}(\log^{2/3}m).
\end{equation}

We sketch a simple proof of (\ref{eq:4.1}). By
(\ref{eq:1.6})
$$\frac{C_{2r}}{C_2}=\prod_{p|r,\,p>2}
\left(1+\frac{1}{p-2}\right).$$
Hence, numbering the primes $p>2$ as
$p_1,\,p_2,\,\cdots$, using the principle of
inclusion-exclusion and letting $[\cdot\cdot]$
denote the integral-part function,
\begin{align} & \qquad
\frac{C_2+\cdots+C_{2m}}{C_2}=m+\sum_j\left[\frac{m}{p_j}
\right]\frac{1}{p_j-2}+\sum_{j,\,k;\,j<k}
\left[\frac{m}{p_jp_k}\right]\cdot\notag \\ &
\qquad\cdot
\left\{\left(1+\frac{1}{p_j-2}\right)\left(1+\frac{1}
{p_k-2}\right)-\frac{1}{p_j-2}-\frac{1}{p_k-2}-1\right\}
+\cdots\notag.
\end{align}
Now simplify, divide by $m$ and let $m\to\infty$. Then by
dominated convergence
\begin{align}
\frac{C_2+\cdots+C_{2m}}{mC_2} &\to
1+\sum\frac{1}{p_j(p_j-2)}+
\sum\frac{1}{p_jp_k(p_j-2)(p_k-2)}+\cdots\notag \\ &
\;\quad =
\prod\left(1+\frac{1}{p_j(p_j-2)}\right)
=\prod\frac{(p_j-1)^2}{p^2_j-2p_j}=\frac{1}{C_2}.\notag
\end{align}
\begin{table} \label{table:3}
\begin{tabular}{lll} 
$q$\hspace{2cm}  & Corresponding $3$-primes $p<500$ &\\
        &                                      & \\
$2,\,4,\,8,\,16$ &
31,\,43,\,109,\,127,\,157,\,223,\,229,\,277,   & \\
     & 283,\,307,\,397,\,433,\,439,\,457,\,499 & \\
$3,\,9,\,24$ & 61,\,67,\,73,\,103,\,151,\,193, & \\   
     & 271,\,307,\,367,\,439,\,499             & \\
$6$  & 7,\,37,\,139,\,163,\,181,\,241,\,307,\,313,& \\
     & 337,\,349,\,379,\,409,\,421,\,439,\,499 & \\
$10$ & 37,\,73,\,79,\,103,\,127,\,139,\,271,   & \\
     & 331,\,349,\,421,\,457,\,463             & \\
$12,\,18$ & 13,\,19,\,79,\,97,\,199,\,211,\,307,& \\
     & 331,\,373,\,439,\,463,\,487,\,499       & \\
$14$  & 13,\,37,\,67,\,79,\,103,\,139,\,157,   & \\
     & 193,\,223,\,379,\,397,\,409,\,439       & \\
$20$ & 7,\,19,\,61,\,97,\,127,\,151,           & \\
     & 193,\,373,\,421,\,457                   & \\ 
$22$ & 7,\,43,\,67,\,73,\,79,\,97,\,103,\,163, & \\
     & 181,\,229,\,331,\,373,\,457             & \\      
     &                                      & 
\end{tabular}
\caption{The `$3$-primes' below $500$}
\end{table}

An elegant proof of (\ref{eq:4.1}) was proposed by
Tenenbaum \cite{Te06}: apply the Wiener--Ikehara theorem
to the Dirichlet series $\sum a(r)r^{-s}$, where $a(r)$
is the multiplicative function $C_{2r}/C_2$. One finds
that the subsequences $\{C_{2hr}\}$ of $\{C_{2r}\}$ have 
mean value $\prod_{p|h,\,p>2}\,p/(p-1)$;
cf.\ Montgomery \cite{Mo71}, Lemma $17.4$. It is plausible
that more generally, the following is true:
\begin{conjecture} \label{con:4.2}
The subsequences of $\{C_{2r}\}$ that correspond to
arithmetic subsequences of the index sequence $\{2r\}$
all have a mean value.
\end{conjecture}
For example, since $\{C_{4r}\}$ has mean value one,
so does the complementary subsequence
$\{C_{4r-2}\}$. The sequence $\{C_{6r}\}$ has
mean value $3/2$, and the sequences $\{C_{6r-2}\}$ and 
$\{C_{6r-4}\}$ should both have mean value $3/4$.

The speculative manuscript \cite{Ko08} suggests an
extension of Proposition \ref{prop:4.1} to the case $k\ge
2$ given by Metatheorem \ref{the:1.2}. Machinery for 
theoretical approach to the metatheorem is developed
in Sections \ref{sec:5}--\ref{sec:8}. 

The constants $C^2_{2r}$ can be computed by using
the Legendre symbol. Table $2$ shows that the average of
these constants with $1\le|r|\le 15$ is about
$0.98$. 

When $k=3$ the computations are more laborious. In order
to obtain a reasonable approximation to
$\ga^3_q$ and $C^3_{2r}$ one has to know the corresponding
$3$-primes (and hence the $0$-primes) up to a suitable
level. Given $q$ we restrict ourselves to primes $p\equiv
1$ (mod $3$) that do not divide $q$. For which $p$ do the
congruences 
\begin{equation} \label{eq:4.3}
n^3\equiv \pm q\;\;(\mbox{mod}\;p)
\end{equation}
have a solution $n\,$? Factorization mod$\;p$ of $n^3-q$
for $|n|\le 100$ will reveal all the $3$-primes $p<200$
and a good many beyond that. To test additional
candidates $p\equiv 1$ (mod $3$) one may use the
following criterion. For given $q$, the
congruences (\ref{eq:4.3}) have a solution (hence three
solutions) if and only if
\begin{equation} \label{eq:4.4}
q^{(p-1)/3}\equiv 1\;\;(\mbox{mod}\;p);
\end{equation}
cf.\ Ireland and Rosen \cite{IR90}, Propositions $7.1.2$
and $9.3.3$. Table $3$ lists the $3$-primes $p<500$ for a
number of values $q$.

For $q=2$ and $q=3$ the constants $\ga^3_q$ were computed
by Bateman and Horn \cite{BH65}, and also by Davenport
and Schinzel \cite{DS66}; the latter constructed
absolutely convergent products. 

For $n^3\pm 4$, $n^3\pm 8$ and $n^3\pm 16$ the $3$-primes
are the same as for $n^3\pm 2$ in Example \ref{exam:3.1}.
Indeed, if $n^3_1\equiv 2$ (mod$\,p$) with $p>2$ and
$n_2\equiv n^2_1$, then $n^3_2\equiv 4$. Conversely, if
$n^3_3\equiv 4$ and $n_4\equiv n^2_3/2$ (mod$\,p$), then
$n^3_4\equiv 2$. 

For $n^3\pm 9$ and $n^3\pm 24$ the $3$-primes are the same
as for $n^3-3$. By the work of Dedekind, the
latter are the primes $p$ for which $4p=a^2+243b^2$. For
$n^3\pm 18$ the $3$-primes are the same as for $n^3\pm
12$. 
\begin{table} \label{table:4}
\begin{tabular}{llll} 
$q=2r$\hspace{1.3cm} & $\ga^3_q$\hspace{1.2cm}
& $C^3_{2r}$    & \\
            &                 &               & \\
$2,\,4,\,16$  & 1.28          & 0.87          & \\
$6$         & 0.82            & 0.98          & \\
$8$         & 0               & 0             & \\ 
$10$        & 1.34            & 1.22         & \\
$12,\,18$   & 0.995           & 1.32          & \\
$14$        & 0.875           & 0.70          & \\
$20$        & 0.73            & 0.58          & \\
$22$        & 0.73            & 0.49          & \\
$24$        & 1.40            & 1.92          & \\
       &                 &               &       
\end{tabular}
\caption{Constants $\ga^3_{2r}$ and $C^3_{2r}$}
\end{table}

Corresponding constants $\ga^3_q$ and $C^3_{2r}$ are given
in Table $4$. Note that $\ga^3_1=0$ and $\ga^3_8=C^3_8=0$
because the corresponding polynomials can be factored.
The average of the constants $C^3_{2r}$ for $1\le r\le
12$ is about $0.93$.

We found $527$ prime pairs $(p,\,p^3+2)$ and
$556$ prime pairs $(p,\,p^3-2)$ with $p<10^5$. With our
imprecise constant $C^3_2\approx 0.87$, the Bateman--Horn
conjecture would give the approximate value
$$(2/3)\cdot 0.87\,{\rm li}_2(10^5)\approx 550.$$

\setcounter{equation}{0} 
\section{Auxiliary functions} \label{sec:5}
For $k\ge 2$ and $r\in{\Bbb Z}\setminus 0$ we
again consider the pair of polynomials
$f_{2r}(n)=\{n,\,n^k+2r\}$. In addition to the counting
function
\begin{equation} \label{eq:5.1}
\pi_{2r}(x) = \pi^k_{2r}(x) = \pi_{f_{2r}}(x)=\#\{p\le
x:\,p^k+2r\;\mbox{prime}\}
\end{equation}
we need the function
\begin{equation} \label{eq:5.2}
\theta_{2r}(x) = \theta^k_{2r}(x) = \sum_{p\le
x;\,p^k+2r\,{\rm prime}}\,\log^2 p.
\end{equation}
Integration by parts will show that for the present
case, the Bateman--Horn Conjecture \ref{con:1.1} 
is equivalent to the asymptotic relation
\begin{equation} \label{eq:5.3}
\theta_{2r}(x)\sim
BH(f_{2r})\,x\quad\mbox{as}\;\;x\to\infty.
\end{equation}
Incidentally, a sieving argument would give
$\theta_{2r}(x)=\cal{O}(x)$; cf.\ Bateman and
Horn \cite{BH65}, Halberstam and Richert
\cite{HR74}, Hindry and Rivoal \cite{HiRi05}.

For the conjecture in the form (\ref{eq:5.3}) we
introduce the Dirichlet series
\begin{equation} \label{eq:5.4}
D_{2r}(s)=D^k_{2r}(s)=\sum_{p;\,p^k+2r\,{\rm
prime}}\,\frac{\log^2
p}{p^{s}}\qquad(s=\si+i\tau,\,\si>1).
\end{equation}
By a two-way Wiener--Ikehara theorem for Dirichlet
series with positive coefficients, relation
(\ref{eq:5.3}) is true if and only if the difference
\begin{equation} \label{eq:5.5}
G_{2r}(s)=D_{2r}(s)-\frac{BH(f_{2r})}{s-1}
\end{equation} 
has `good' boundary behavior as $\si\searrow 1$.
That is, $G_{2r}(\si+i\tau)$ should tend to a
distribution
$G_{2r}(1+i\tau)$ which is locally equal to a
pseudofunction. By a pseudofunction we mean the
distributional Fourier transform of a bounded function
which tends to zero at infinity; it
cannot have poles. A pseudofunction may be characterized
as a tempered distribution which is locally given by
Fourier series whose coefficients tend to zero; see
\cite{Ko05}. In particular $D_{2r}(s)$ itself would have
to show pole-type behavior, with residue $BH(f_{2r})$, for
angular approach of $s$ to $1$ from the right;
there should be no other poles on the line $\{\si=1\}$.

In the following we have to use repeated
complex integrals related to those in \cite{Ko07}.

\setcounter{equation}{0}
\section{Complex integral for a sieving function}
\label{sec:6} 
Our integrals involve sufficiently smooth even
sieving  functions $E^\la(\nu)=E(\nu/\la)$ depending
on a parameter $\la>0$. The basic functions $E(\nu)$
have $E(0)=1$ and support $[-1,1]$; it is required
that $E$, $E'$ and $E''$ be absolutely continuous,
with $E'''$ of bounded variation. One may for example 
take
\begin{align} \label{eq:6.1}
E^\la(\nu) &= 
\frac{3}{4\pi}\int_0^\infty\frac{\sin^4(\la t/4)}
{\la^3(t/4)^4}\cos \nu t\,dt
\notag \\ &=\left\{\begin{array}{ll}
1-6(\nu/\la)^2+6(|\nu|/\la)^3 & \mbox{for
$|\nu|\le\la/2$},\\ 2(1-|\nu|/\la)^3 &
\mbox{for
$\la/2\le|\nu|\le\la$},\\ 0 & \mbox{for
$|\nu|\ge\la$.}
\end{array}\right.
\end{align}
An important role is played by a Mellin transform
associated with the Fourier transform of the kernel
$E^\la(\nu)=E(\nu/\la)$. For $0<x={\rm Re}\,z<1$
\begin{align} \label{eq:6.2}
M^\la(z) & \stackrel{\mathrm{def}}{=}
\frac{1}{\pi}\int_0^\infty \hat E^\la(t)t^{-z}dt
= \frac{2}{\pi}\int_0^\infty t^{-z}dt\int_0^\la
E^\la(\nu)(\cos t\nu)d\nu
\notag \\
&= \frac{2}{\pi}\int_0^\la E(\nu/\la)d\nu
\int_0^{\infty-}(\cos\nu t)t^{-z}dt\notag\\ &=
\frac{2}{\pi}\Ga(1-z)\sin(\pi z/2)\int_0^\la
E(\nu/\la)\nu^{z-1}d\nu \\ &=
\frac{2\la^z}{\pi}\Ga(1-z)\sin(\pi z/2)\int_0^1
E(\nu)\nu^{z-1}d\nu \notag \\ &=
\frac{2\la^z}{\pi}\Ga(-z-3)\sin(\pi z/2)\int_0^{1+}
\nu^{z+3}dE'''(\nu).\notag
\end{align}
The Mellin transform extends to a
meromorphic function for $x>-3$ with simple poles at 
the points $z=1,\,3,\,\cdots$. The residue of the pole
at $z=1$ is $-2(\la/\pi)A^E$ with $A^E=\int_0^1
E(\nu)d\nu$, and $M^\la(0)=1$. Setting $z=x+iy$ (and
later $w=u+iv$), the standard order estimates
\begin{equation} \label{eq:6.3}
\Ga(z)\ll |y|^{x-1/2}e^{-\pi|y|/2},\quad\sin(\pi z/2)
\ll e^{\pi|y|/2}
\end{equation}
for $|x|\le C$ and $|y|\ge 1$ imply the useful
majorization
\begin{equation} \label{eq:6.4}
M^\la(z)\ll \la^x (|y|+1)^{-x-7/2}\quad\mbox{for}\;\;
-3<x\le C,\;|y|\ge 1.
\end{equation}
\noindent
{\scshape Repeated complex integral for
$E^\la(\al-\be)$}. We write $L(c)$ for the `vertical
line' $\{x=c\}$; the factor $1/(2\pi i)$ in complex
integrals will be omitted. Thus
$$\int_{L(c)}f(z)dz\stackrel{\mathrm{def}}{=}
\frac{1}{2\pi i}\int_{c-i\infty}^{c+i\infty}f(z)dz.$$
Since it is important for us to have absolutely 
convergent integrals, we often have to replace a line
$L(c)$ by a path $L(c,B)=L(c_1,c_2,B)$ with suitable
$c_1<c_2$ and $B>0$:
\vspace{-5mm}
\begin{figure}[htb]
$$\includegraphics[width=4cm]{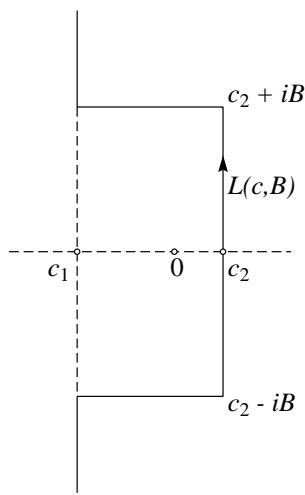}$$
\vspace{-5mm}
\caption{The path $L(c_1,c_2,B)$ \label{fig:1}}
\end{figure} 
\begin{equation} \label{eq:6.5}
L(c,B)=\left\{\begin{array}{lllll}
\mbox{$\quad$the half-line} & \mbox{$\{x=c_1,\,-\infty
<y\le -B\}$}\\
\mbox{$+\;$the segment} & \mbox{$\{c_1\le x\le
c_2,\,y=-B\}$}\\
\mbox{$+\;$the segment} & \mbox{$\{x=c_2,\,-B\le y\le 
B\}$}\\ \mbox{$+\;$the segment} & \mbox{$\{c_2\ge x\ge
c_1,\,y=B\}$}\\
\mbox{$+\;$the half-line} & \mbox{$\{x=c_1,\,B\le
y<\infty\}$;}  
\end{array}\right.
\end{equation}
cf.\ Figure \ref{fig:1}. Thus, for example, 
$$\cos\al = \int_{L(c,B)}\Ga(z)\al^{-z}\cos(\pi
z/2)dz \qquad(\al>0),$$
with absolute convergence if $c_1<-1/2$ and $c_2>0$.
Similarly for $\sin\al$. For the combination
$$\cos(\al-\be)t=\cos\al t\cos\be t+\sin\al
t\sin\be t$$  
with $\al,\,\be,\,t>0$, one can now write down an
absolutely convergent repeated integral. 
In \cite{Ko07} it was combined with (\ref{eq:6.2}) to
obtain a repeated complex integral for the sieving
function $E^\la(\al-\be)$ in which $\al$ and
$\be$ occur separately. Taking
$-3<c_1+c'_1<0$, $c_2,\,c'_2>0$ and $c_2+c'_2<1$ one has
\begin{align} \label{eq:6.6}
E^\la(\al-\be) & = \int_{L(c,B)}\Ga(z)\al^{-z}dz
\int_{L(c',B)}\Ga(w)\be^{-w}\,\cdot\notag
\\ & \quad\cdot M^\la(z+w)\cos\{\pi(z-w)/2\}\,dw.
\end{align}
To verify the absolute convergence of the repeated
integral one would substitute $z=x+iy$,
$w=u+iv$ and then use the inequalities (\ref{eq:6.3}),
(\ref{eq:6.4}) and $\cos\{\pi(z-w)/2\}\ll
e^{\pi(|y|+|v|)/2}$, together with a simple lemma:
\begin{lemma} \label{lem:6.1}
For real constants $a,\,b,\,c$, the function
$$\phi(y,v)=(|y|+1)^{-a}(|v|+1)^{-b}(|y+v|+1)^{-c}$$
is integrable over ${\Bbb R}^2$ if and only if $a+b>1$,  
$a+c>1$, $b+c>1$ and $a+b+c>2$. 
\end{lemma}

\setcounter{equation}{0} 
\section{Repeated integral for a function
$T^\la_k(s)$} \label{sec:7}
Taking $k\ge 2$ and using paths specified below, we
consider integrals
\begin{align} \label{eq:7.1}
T^\la_k(s) & = \int_{L(c,B)}\Ga(z-s)
\frac{\ze'(kz)}{\ze(kz)}\,dz
\int_{L(c',B)}\Ga(w-s)\frac{\ze'(w)}{\ze(w)}\,
\cdot \notag \\ & \quad\; \cdot
\,M^\la(z+w-2s)\cos\{\pi(z-w)/2\}dw.
\end{align}
The case $k=1$ was used in \cite{Ko07} to study
prime pairs $(p,\,p+2r)$. Proceeding in a
similar way, we use the Dirichlet series 
$$\frac{\ze'(Z)}{\ze(Z)}=-\sum
\frac{\La(n)}{n^Z}=-\sum(\log p)
\bigg(\frac{1}{p^Z}+\frac{1}{p^{2Z}}+
\cdots\bigg)$$
and formula (\ref{eq:6.6}) to obtain the (as yet formal)
expansion
\begin{align} 
T^\la_k(s) &=
\sum_{h,\,j}\,\La(h)\La(j)h^{-ks}j^{-s} E^\la(h^k-j)
\notag \\ &= 
\sum_{0\le|d|\le\la}\,\sum_h\,\La(h)\La(h^k-d)h^{-ks}
(h^k-d)^{-s}E^\la(d) \notag \\ &=
\sum_{0\le
|2r|\le\la}\,\sum_h\,\La(h)\La(h^k-2r)h^{-2ks}
E^\la(2r) + H^\la_1(ks),\notag
\end{align}
where $H^\la_1(Z)$ is holomorphic for $X>0$. Indeed,
$(h^k-d)^{-s}$ may be approximated by
$h^{-ks}$, and for odd numbers $d$, the product
$\La(h)\La(h^k-d)$ can be $\ne 0$ only if $h=2$ or
$h^k-d=2$. The expansion can be used to define
$T^\la_k(s)$ as a holomorphic function for
$\si>1/(2k)$. 

Having even $d=2r$, the principal
contributions to the expansion come from the cases
where $h$ is a prime $p$, and either $h^k-2r$ is the prime
power $p^k$ (if $r=0$), or a prime $q$ (if $r\ne 0$).
In the latter case $\La(h^k-2r)=\log q=\log(p^k-2r)$
is well-approximated by $k\log p$. Using
(\ref{eq:5.4}) one thus finds that
\begin{equation} \label{eq:7.2}
T^\la_k(s) = D_0(2ks) + k\sum_{0<|2r|\le
\la}\,E(2r/\la)D_{2r}(2ks)+H^\la_2(ks),
\end{equation}
where $D_0(Z)=\sum_p\,(\log^2 p)/p^Z$ and
$H^\la_2(Z)$ is holomorphic for $X>1/3$. Comparison of
$D_0(Z)$ with $(d/dz)\sum\La(n)/n^Z$  shows that
$D_0(Z)$ is holomorphic for $X>1/4$, except for
purely quadratic poles at $Z=1$, $Z=1/2$ and the
complex zeros $Z=\rho$ of $\ze(Z)$. 

Assuming Riemann's Hypothesis (RH) for simplicity, one may
in (\ref{eq:7.1}) take $c_1=(1/4)+\eta$,
$c_2=(1/2)+\eta$ and $c'_1=(1/2)+\eta$, $c'_2=1+\eta$
with small $\eta>0$. Varying $\eta$, the integral will
represent $T^\la_k(s)$ as a holomorphic function for
$3/8<\si<1/2$ and $|\tau|<B$. Indeed, for given
$s$ and small $\eta$ there will be no singular points on
the paths. Absolute convergence (locally uniform in $s$)
may be verified with the aid of Lemma \ref{lem:6.1}, using
(\ref{eq:6.3}), (\ref{eq:6.4}) and the fact that
$(\ze'/\ze)(Z)$ grows at most logarithmically in $Y$ on
vertical lines $\{X=d\}$ with $d\ne 1/2$; cf.\ Titchmarsh
\cite{Ti86}. Thus on the remote parts of the paths, the
integrand is majorized by
$$C(\la)|y|^{c_1-\si-1/2}(\log|y|)|v|^{c'_1-\si-1/2}
(\log|v|)(|y+v|+1)^{-c_1-c'_1+2\si-7/2}.$$
Similar estimates will enable us to move the paths of
integration; cf.\ \cite{Ko07}.

Starting with $s$ and the paths in (\ref{eq:7.1}) as
above, we now move the $w$-path $L(c',B)$ across the poles
at the points $w=1$, $\rho$ and $s$ to the path $L(d,B)$
with $d_1=-1/2$ and $d_2=0$. Then the residue theorem
gives
\begin{equation} \label{eq:7.3}
T^\la_k(s) =
\int_{L(c,B)}\cdots\,dz\int_{L(d,B)}\cdots\,dw
+U^\la_k(s)=T^{\la *}_k(s)+U^\la_k(s),
\end{equation}
say, where 
\begin{equation} \label{eq:7.4}
U^\la_k(s) = \int_{L(c,B)}
\Ga(z-s)\frac{\ze'(kz)}{\ze(kz)}\,J(z,s)dz,
\end{equation}
with 
\begin{align} \label{eq:7.5}
J(z,s) &=
-\Ga(1-s)M^\la(z+1-2s)\cos\{\pi(z-1)/2\}
\notag \\ & \quad  
+\sum_\rho\,\Ga(\rho-s)M^\la(z+\rho-2s)
\cos\{\pi(z-\rho)/2\} \\ & \quad
+\frac{\ze'(s)}{\ze(s)}\,M^\la(z-s)\cos\{\pi(z-s)/2\}.
\notag
\end{align}

Observe that the apparent poles of $J(z,s)$ at the points
$s=1$ and $s=\rho$ cancel out.
For given $s$ with $3/8<\si<1/2$ and $|\tau|<B$, and
for suitably small $\eta$, the function $J(z,s)$ is
holomorphic in $z$ on and between the paths $L(c,B)$ and
$L(d,B)$. We now move the path $L(c,B)$ in (\ref{eq:7.4})
to $L(d,B)$. Picking up residues at $z=s$, $1/k$ and the
zeros $\rho'/k$ of $\ze(kz)$, one finds that
\begin{align} \label{eq:7.6} 
U^\la_k(s) &=
\int_{L(d,B)}\,\Ga(z-s)\frac{\ze'(kz)}{\ze(kz)}\,
J(z,s)dz + V^\la_k(s) \notag \\ &=
U^{\la *}_k(s)+V^\la_k(s),
\end{align}
say, where 
\begin{align} \label{eq:7.7}
V^\la_k(s) &= \frac{\ze'(ks)}{\ze(ks)}\,J(s,s) 
-(1/k)\Ga\{(1/k)-s\}J(1/k,s)\notag \\ &
\quad + \sum_{\rho'}\,(1/k)\Ga\{(\rho'/k)-s\}J(\rho'/k,s).
\end{align}

The single integral for $U^{\la *}_k(s)$ in (\ref{eq:7.6})
defines a holomorphic function for $0<\si<1/2$. Varying
$B$, the same is true for the repeated integral defining
$T^{\la *}_k(s)$ in (\ref{eq:7.3}). For 
verification one may use Lemma \ref{lem:6.1} or an analog
for the integral of a sum; the number of points $\rho$
with $n-1<{\rm Im}\,\rho\le n$ is $\cal{O}(\log n)$.

We know that $T^\la_k(s)$ is holomorphic in the  strip
$$\cal{S}=\{1/(2k)<\si<1/2\}.$$
For our work we have to know the {\it behavior} of
$T^\la_k(s)$ {\it near the boundary line}
$\{\si=1/(2k)\}$. This is determined by the sum
$V^\la_k(s)$ in (\ref{eq:7.7}), which in view of
(\ref{eq:7.5}) splits into nine separate terms. All but
one of these clearly represent meromorphic
functions in the closed strip $\overline{\cal{S}}$. The
exception is the function defined by the double series
arising from the third term in (\ref{eq:7.7}):
\begin{align} \label{eq:7.8}
\Si^\la_k(s) &\stackrel{\mathrm{def}}{=}
\sum_{\rho,\,\rho'}\,(1/k)\Ga(\rho-s)\Ga\{(\rho'/k)-s\}
\,\cdot \notag \\ & \qquad
\cdot M^\la(\rho-2s+\rho'/k)\cos\{\pi(\rho-\rho'/k)/2\}.
\end{align}
By a discrete analog of Lemma \ref{lem:6.1} the
double series is absolutely convergent for
$(1+1/k)/4<\si<1/2$. The analysis below will show that
the sum has an analytic continuation [also denoted
$\Si^\la_k(s)$] to $\cal{S}$; see (\ref{eq:8.4}).

\setcounter{equation}{0} 
\section{Behavior of
$T^\la_k(s)$ near the line $\{\si=1/(2k)\}$}
\label{sec:8}
We start with the second term of $V^\la_k(s)$ in
(\ref{eq:7.7}). By (\ref{eq:7.5}) the factor $J(1/k,s)$
is holomorphic in $\overline{\cal{S}}$, except for a
simple pole at $s=1/(2k)$ due to the pole of $M^\la(Z)$
for $Z=1$ with residue $-2(\la/\pi)A^E$; cf.\
(\ref{eq:6.2}). Also taking into account the other factor
$-(1/k)\Ga\{(1/k)-s\}$, a short calculation gives the
principal part of the pole at
$s=1/(2k)$ as
\begin{equation} \label{eq:8.1}
\frac{(1/k)A^E\la}{s-1/(2k)},\quad\mbox{where}\;\;
A^E=\int_0^1 E(\nu)d\nu.
\end{equation} 
There is also a pole at $s=1/k$, but it is cancelled by a
pole of the first term in $V^\la_k(s)$. That term involves
$J(s,s)$, which by (\ref{eq:7.5}) is holomorphic in
$\overline{\cal{S}}$, and $(\ze'/\ze)(ks)$, which
besides $s=1/k$ has poles at the points $\rho'/k$ on the
line $\si=1/(2k)$. The latter have principal parts
\begin{equation} \label{eq:8.2}
\frac{(1/k)J(\rho'/k,\rho'/k)}{s-\rho'/k}.
\end{equation}

The third term of $V^\la_k(s)$ involves an infinite
series of products. The factors $J(\rho'/k,s)$ are
holomorphic in $\overline{\cal{S}}$, but the factors
$(1/k)\Ga\{(\rho'/k)-s\}$ introduce poles at the points
$s=\rho'/k$. The poles in the products have
principal parts 
\begin{equation} \label{eq:8.3}
\frac{-(1/k)J(\rho'/k,\rho'/k)}{s-\rho'/k},
\end{equation}
hence these poles cancel those given by (\ref{eq:8.2}).
The final term of $J(\rho'/k,s)$ leads to the
function $\Si^\la_k(s)$ defined by the double series in
(\ref{eq:7.8}).
\begin{summary} \label{sum:8.1}
Assume RH. Combination of (\ref{eq:7.2}) and the
subsequent results in Sections \ref{sec:7} and
\ref{sec:8} shows that in the strip
$\cal{S}=\{1/(2k)<\si<1/2\}$,
\begin{align} \label{eq:8.4}
T^\la_k(s) &=
D_0(2ks)+k\sum_{0<|2r|\le\la}\,E(2r/\la)D_{2r}(2ks)
  + H^\la_2(ks) \notag \\ &=
\frac{(1/k)A^E\la}{s-1/(2k)}
+\Si^\la_k(s)+H^\la_3(s),
\end{align}
where $H^\la_2(ks)$ and $H^\la_3(s)$ are
holomorphic for $1/(2k)\le\si<1/2$.
\end{summary} 
We now focus on the difference $\Si^\la_k(s)-D_0(2ks)$,
which by (\ref{eq:8.4}) can be considered as a
holomorphic function in $\cal{S}$. How does it behave
as $s$ approaches the line $\{\si=1/(2k)\}\,$? The
function $D_0(2ks)$ has a purely quadratic pole at
$s=1/(2k)$; see (\ref{eq:7.2}). By sieving, the functions
$D_{2r}(2ks)$ cannot have a pole at $s=1/(2k)$ of higher
order than the first; cf.\ Section \ref{sec:5}, hence
$\Si^\la_k(s)$ must cancel the quadratic pole of
$D_0(2ks)$. On the basis of the Bateman--Horn
conjecture in (\ref{eq:5.3}) it is {\it plausible} that
the functions $D_{2r}(2ks)$ do have first-order poles at
$s=1/(2k)$, with respective residues $BH(f_{2r})/(2k)$;
cf.\ (\ref{eq:5.5}). 

Assuming (\ref{eq:5.3}), what can we say about the
residue of $\Si^\la_k(s)-D_0(2ks)$ for $s\searrow
1/(2k)\,$? By (\ref{eq:8.4}) and (\ref{eq:8.1}) it will
be equal to
\begin{equation} \label{eq:8.5}
R_k(\la)=k\sum_{0<|2r|\le\la}\,E(2r/\la)BH(f_{2r})/(2k)-
(\la/k)\int_0^1 E(\nu)d\nu. 
\end{equation}
Now it is plausible that this residue is $o(\la)$ as
$\la\to \infty$. Indeed, $\la$ occurs in the terms of
$\Si^\la_k(s)$ only as a factor
$\la^{\rho-2s+\rho'/k}$. Cf.\ the case of
$T^\la_1(s)$ and the sum $\Si^\la_1(s)$ in
\cite{Ko07}, where one dealt with
ordinary prime pairs $(p,\,p+2r)$, so that
$BH(f_{2r})=2C_{2r}$ and it {\it is known} that
$R_1(\la)=o(\la)$; see (\ref{eq:4.1}). By analogy
assuming $R_k(\la)=o(\la)$, and letting
$E(\nu)\le 1$ approach the constant function
$1$ on $[0,1]$, it follows from (\ref{eq:8.5}) that 
\begin{equation} \label{eq:8.6}
\sum_{0<|r|\le\la/2}\,BH(f_{2r})/2\sim\la/k
\quad\mbox{as}\;\;\la\to\infty.
\end{equation}
Hence by (\ref{eq:1.3}) or (\ref{eq:1.9}), the numbers
$C^k_{2r}=(k/2)BH(f_{2r})$ should have mean value
$1$, as asserted in Metatheorem \ref{the:1.2}. 
\begin{remark} \label{rem:8.2}
By more refined treatment of $T^\la_k(s)$ the
conclusion can be obtained without RH; cf.\ the analysis
of $T^\la_1(s)$ in \cite{Ko07}.
\end{remark}

\setcounter{equation}{0} 
\section{The Bateman--Horn constants $\ga^k_q$}
\label{sec:9} 
Heuristics based on the relevant Bateman--Horn
conjectures suggest that the constants
$\ga^k_q$ also have mean value one. This is supported by
numerical evidence; a rough computation gives the average
of $\ga^2_q$ for $1\le q\le 20$ as $0.99$ and for
$-20\le q\le -1$ as $1.03$. There are corresponding
results for even $q$; cf.\ Tables $2$ and $4$.

For the study of $\ga^k_q$ one may introduce a related
Dirichlet series. By (\ref{eq:1.11}),
$$\ga^k_q = \lim_{s\searrow 1}\,
\prod_p\,\bigg(1-\frac{1}{p^s}\bigg)^{-1}\prod_p\,
\bigg(1-\frac{\nu^k_q(p)}{p^s}\bigg)
= \lim_{s\searrow 1}\,\ze(s)G^k_q(s),$$
say. Kurokawa \cite{Ku86a}, cf.\ \cite{Ku86b}, has studied
the general product
$$Z(s,f)=\ze^m(s)\prod_p\,\{1-N_f(p)p^{-s}\},$$
which is related to the product for $BH(f)$ in 
(\ref{eq:1.3}). See also Conrad \cite{Con03}.

The mean value one of prime $n$-tuple constants plays a
role in recent work of Kowalski \cite{Kow08}.

\bigskip

\noindent{\scshape KdV Institute of Mathematics,
University of Amsterdam, \\
Plantage Muidergracht 24, 1018 TV Amsterdam,
Netherlands}

\noindent{\it E-mail}: {\tt fjvdbult@science.uva.nl},
{\tt korevaar@science.uva.nl}

\enddocument